\documentstyle[12pt]{article}

\begin{document}

\author{Mai Gehrke and David Pengelley \\
Mathematical Sciences\\
New Mexico State University\\
Las Cruces, NM 88003\\
mgehrke@nmsu.edu davidp@nmsu.edu\\
Preprint: to appear in {\em Calculus: The Dynamics of Change}, pp. 20--23, \\
MAA Notes {\bf 39}, Mathematical Association of America, 1996}
\title{Towards active processes for teaching and learning}
\date{}
\maketitle

What are the real goals of `reforming' calculus teaching? Some commonly
accepted aims are making calculus more relevant and understandable to
students, making it `lean and lively', and having students acquire and meld
the tools of calculus to solve multi-step or open-ended problems. However,
all these aims are actually just part of an overarching goal: having
students actively involved and taking initiative in their own learning, in
fact learning how to learn for themselves. Our students should take
responsibility for and charge of their own learning, developing their own
process for becoming independent learners, and thus end up with a sense of
personal ownership of the results of their labor. The commonly accepted aims
listed above are good first steps on this road towards helping students
become active learners: by creating material relevant to students' curiosity
we make it possible and attractive for them to take interest; by making the
syllabus lean enough, we allow time for the high level of absorption
inherent in active learning; by having students solve larger contextual
problems rather than template problems floating in a vacuum, calculus
becomes understandable and useful in a more real way.

So what are the means for achieving the broad goal of making students active
learners? Is it accomplished by incorporating technology, or by adopting a
`reformed' text, or through substantial individual or group projects? All
these can be useful pedagogical tools in trying to reach the broad goal.
However, if we lose sight of the overarching goal, assuming that adopting
one or more of these specific tools is the essence of reform, then we will
fail; when the tools themselves become the focus, we depart from the
necessary primary emphasis on the active involvement process. Another reason
not to restrict one's view to a specific tool is that if our students are
truly becoming active learners, due to their individuality they will each
develop personal learning processes which thrive on different tools.

Since our ultimate goal for students is a focus on their learning process,
the tools we use should always be treated merely as objects within that
process. Their importance and usefulness should be kept in perspective, and
thus the teacher's focus should also shift away from objects. Objects such
as text materials, computers, labs, lectures, reading, writing, projects,
group work, exams, and portfolios are colors on the palette from which
teacher and students can choose during the ongoing process of each student's
learning. Ultimately, we believe this shift will need to encompass all
levels of the educational structure: students, teachers, and the educational
community as a whole.

Our professional community has the same responsibility to teachers that
teachers have to students: to create an environment in which teachers will
naturally evolve an ongoing active and creative teaching process. In
particular, teachers must be encouraged to show individual initiative in
order for reform to succeed. Each teacher must `make it their own',
selecting from a community palette of ideas and resources.

Our own experience is that in a department where collaboration in teaching
innovation is strongly encouraged, while individual instructors still make
their own decisions about how to teach their classes, we have an unusually
high level of participation in `reform' without forcing involvement. Our
department values the varied contributions of many teachers, and these are
continually being synthesized into our own current version of reform;
individuals and groups of teachers communicate their fresh approaches, often
resulting in the incorporation of their ideas into the reform undertaken by
others. What our department `lacks' is an orthodoxy of reformed teaching
materials or other objects, since we recognize that individual teachers will
go about creating process differently.

Giving teachers the freedom to develop their own personalized reform is what
incites initiative from them, and thus ultimately from their students.
Faculty can then become active teachers, and the multitude of approaches
they develop will naturally induce healthy cross-fertilization. Of course
change in a given department may not start with individual initiative from
each teacher, but can nevertheless evolve into an organic environment for
reform provided the focus of change is not so rigidly tied to certain
objects of reform that it frustrates the development of individual process
for teachers.

Thus the community should nurture a fluid atmosphere, in which adopting
pre-existing reform materials can go hand in hand with initiating an
individual teacher's process of change, and should provide a library of
resources from which teachers can create personalized pedagogical tools for
their own teaching. With this kind of individual growth, change will happen
more slowly than with a superficial imposition or adoption of prepackaged
reform, but will surely be more longlasting, and more faithful to the goals
above; a teacher creating her/his own materials, or modifying those created
by others, will be an active teacher, whereas simply using prepackaged
materials will not stimulate a teacher to emerge from passivity.

What we have experienced at New Mexico State University gives an example of
how this can happen on a sizeable scale. We began with a small group of
faculty initiating change. Then a larger group of faculty found themselves
enticed into getting involved, and this subsequently shifted the pedagogical
nature of our efforts. In this way the group of faculty involved has kept
growing, and our direction has evolved with this growth. Thus our `reform'
has happened in stages, each reflecting a new horizon which only became
visible at the end of a previous stage. Our evolution can provide an example
of the dynamics of individual and departmental cultural change.

For us, change began in 1987 as a response to unsatisfactory student
performance in calculus. Two faculty members came up with some basic ideas
for improving the situation. One was to grab students' attention --- we were
forever hearing from students that they could not work on their mathematics
homework because they had assignments due in other classes. A second idea
was to have students do some real mathematics --- problems that they would
solve and explain as opposed to ones to which they merely supplied an
answer. Basically, we wanted them to think and we needed a way to encourage
them to do that.

Projects were designed to get students to think for themselves on major
multistep, take-home problems, working individually or in groups. We hoped
to alter fundamentally students' view of what mathematics is all about and
simultaneously build their self-confidence in what they could achieve
through imaginative, theoretical thinking. The projects resemble
mini-research problems. Most of them require creative thought and all of
them engage students' analytic and intuitive faculties, often weaving
together ideas from many parts of calculus. While many of the projects are
couched in seemingly real-world settings, often with engaging story lines,
they are all in a sense theoretical. One cannot do them without an
appreciation of the ideas behind the method. Students must decide what the
problem is about, what tools from the calculus they will use to solve it,
find a strategy for its solution, and present their findings in a written
report. This approach yields an amazing level of sincere questioning,
energetic research, dogged persistence, and conscientious communication from
students. Moreover, our own opinions of our students' capabilities
skyrocketed as they rose to the challenges presented by these projects, and
some other faculty and graduate teaching assistants were smitten and wanted
to get involved.

Even though the idea of having students work on projects seemed a
revolutionary idea at its conception, it was a small enough step that a
number of faculty felt comfortable about incorporating a project or two in
their courses. The new teachers wanted to create their own projects, or
modify old ones, each bringing a unique perspective to what a `project'
should be, and thus they became active in reform. The motivations and types
of projects written by this conglomeration of people varied and added
breadth and scope to the nature and efficacy of using projects in teaching.
Over 100 projects were developed by five faculty during this period, and
published in the MAA\ book {\em Student Research Projects in Calculus},
along with several chapters detailing the logistics of assigning projects
and advice for instructors.

In 1990 the program expanded and branched in various directions. Numerous
other faculty in the department volunteered to use calculus projects in
their classes, and we began the development of a discovery-project based
vector calculus and differential equations curriculum (in which a continuous
sequence of discovery projects forms the context for learning all the
material of the course); we also started a collaborative program with local
high school teachers to bring projects into high school mathematics courses.

As new faculty became involved in teaching with projects, they injected
fresh ideas into the program and the projects approach itself evolved.
Although introducing projects was a valid first step, we realized this had
created somewhat schizophrenic courses in which students worked on projects
outside class, while the classroom continued to function in a traditional
style. Even though we felt that the activities involved in working on
projects were effective in stimulating students to think and to learn
mathematics, our day-to-day classroom activities remained largely unchanged.
This provoked a new stage in our development.

Thus in 1991 a group of faculty pioneered a major new emphasis on
cooperative self-learning both in and out of the classroom, developing
structured in-class assignments called `themes'. A distinct change is that
themes are used to introduce the core material of the course and much class
time is spent working on them, with less time on lecture, whereas the
projects were completed outside class and contained material over and above
day-to-day course work.

In a theme assignment, students learn and write about core course material
while working in groups with the instructor serving as a resource. When
themes were first assigned, students completed a written theme report every
week. Experience has tempered this pace somewhat, and we are now assigning
three to six themes per semester. Today, several instructors are blending
the theme approach with the discovery-project methods developed in vector
calculus and differential equations courses. Other ideas, such as class
discussions, student presentations, and mastery skills exams, are being
tried also. These somewhat independent directions seem to be
cross-fertilizing each other's growth.

In retrospect, we see that each of the tools we develop leads to new
pedagogical challenges. For instance, we were pleased with the high level of
student initiative and achievement that projects elicited, but we wanted to
get away from the passive role of our students during a lecture. With themes
and writing assignments, students were active in the classroom, but we
realized the teacher should be more than just a resource for individual
students or small groups to call on; in fact, this placed the teacher in too
passive a role. The teacher should be providing leadership to the class as a
whole, in order to take advantage of having all the students and teacher
together. Sometimes while students were working on themes, most groups would
generate a common question, which naturally led to a whole-class discussion
moderated and guided by the instructor. In fact the guided class discussions
based on students' questions arising from their active work emerged as one
of the most successful and productive aspects of this student-centered
classroom. We now view such guided class discussions as an important tool in
their own right, and we have found other student activities which benefit
from and enable these discussions. For instance, another of our aims has
been for our students to become capable and active readers of mathematics.
This requires breaking the vicious cycle in which instructors lecture text
material to students because they know students don't actually learn it from
reading, and students have little incentive to read because they know their
instructor will lecture it to them. We have found that if we demand students
read in advance, and write commentary and questions about their reading,
then these questions can form the basis for active class discussion,
bypassing the vicious cycle and leading to more productive and satisfying
classroom learning.

Theme assignments have also prompted us to incorporate structured means of
improving student skills at mathematical and prose report writing. After
incorporating handouts on writing, and learning how to guide students in
honing their writing skills, we have seen an incredible improvement in their
ability to write. Reading and writing in mathematics have emerged as
important features of reform at both the undergraduate and graduate level,
and can be viewed as a new stage which has spread far beyond our calculus
courses. These new emphases have merged with innovative efforts of other
faculty who were never even involved in our calculus reform program.

Of course our means and methods for grading have also changed drastically
over this period of reform. When we introduced projects and themes as
learning tools, we also used them as an important means of evaluation. In
comparison we found that traditional exams have little to do with learning,
and we now primarily value means of evaluation that are also learning tools.
Our methods of grading began to change as well, since already the projects
required us to learn how to evaluate written reports and group work. We came
to realize that detailed numerical grading was poorly suited to grading
large written reports, so we have been learning how to evaluate student work
in a more holistic fashion. A benefit of qualitative holistic grading is
that students get specific feedback on how to improve their written work. At
first we worried that students would feel uncomfortable not having points
attached to every aspect of their work, but we found that they readily
accept and appreciate qualitative feedback and evaluation; it is easier for
them to see the qualitative nature of the distinction between A and B work,
provided the criteria are clearly explained, than to understand the
difference between grades of 89 and 90.

After seeing the benefits of holistic evaluation for individual assignments,
it is natural to consider extending this methodology to evaluation of the
totality of a student's work. Recently this has led some of us to a
portfolio approach, in which the student prepares a showcase of their entire
work for the semester, and this portfolio is evaluated as a whole at the
end. In contrast, traditional grading is by nature fragmented, encouraging a
disconnected view in the student of both the course and their own work. By
assuming responsibility for collecting, organizing, and presenting all their
course work in a portfolio, students become aware of the big picture in both
the subject matter and their own performance.

Looking back on all these changes, we see that they involved a relinquishing
of total control. While this can be a frightening prospect, it is necessary
if students are to assume more responsibility and control of their own
learning. Fortunately, if the balance of control is shifted gradually from
teacher to students, through a slow process of evolution, total loss of
control may be avoided. The reward is the opening of new vistas for both
teaching and learning, in which the instructor becomes an expert guide,
facilitator, and coordinator. Even though the original purpose is to improve
student learning, there is a tremendous revitalizing benefit for teachers,
as our interaction with students and colleagues becomes more rewarding, and
the results of our efforts become more meaningful.

Each of the specific changes and pedagogical discoveries that we have made
along the way has in a sense forced itself upon us as an inescapable
outgrowth of a previous change. This process, and the collegial atmosphere
that has made it possible, are in our minds the essential features of reform
as we have experienced it. Alan Schoenfeld, in the preface to the recently
published book {\em Mathematical Thinking and Problem Solving}, referred to
our initial seed, namely student projects, as a Trojan Mouse, and that is
truly what it has been, subversively driving the scope of change far beyond
what we could originally imagine. The other essential feature of our reform
is the atmosphere of faculty collaboration in teaching innovation; it has
nurtured almost everything we have accomplished, and in a way which has
fostered individual faculty ownership of both the process and the results,
anchoring it deeply in the fabric of our department community.

\end{document}